\documentclass[11pt]{article}
\usepackage{amsmath}
\usepackage{amsthm}
\usepackage{amsfonts}
\usepackage{amssymb}

\makeatletter
\def\labelenumi{(\makebox[8pt]{\@arabic\c@enumi})}
\def\theenumi{\@arabic\c@enumi}
\makeatother
\theoremstyle{plain}
\newtheorem{thm}{Theorem}
\newtheorem{lem}{Lemma}

\theoremstyle{definition}
\newtheorem{remark}{Remark}
\numberwithin{equation}{section}
\def\disc{\mathop{\rm disc}}

\def\tr{\mathop{\rm tr}}

\def\sgn{\mathop{\rm sgn}}

\def\Sym{\mathop{\rm Sym}\nolimits}

\makeatletter
\def\multilimits@{\bgroup
  \Let@
  \restore@math@cr
  \default@tag
 \baselineskip\fontdimen10 \scriptfont\tw@
 \advance\baselineskip\fontdimen12 \scriptfont\tw@
 \lineskip\thr@@\fontdimen8 \scriptfont\thr@@
 \lineskiplimit\lineskip
 \vbox\bgroup\ialign\bgroup\hfil$\m@th\scriptstyle{##}$\hfil\crcr}
\def\Sb{_\multilimits@}
\def\Sp{^\multilimits@}
\def\endSb{\crcr\egroup\egroup\egroup}

\def\smallmatrix{\null\,\vcenter\bgroup
 \Let@\restore@math@cr\default@tag
 \baselineskip6\ex@ \lineskip1.5\ex@ \lineskiplimit\lineskip
 \ialign\bgroup\hfil$\m@th\scriptstyle{##}$\hfil&&\thickspace\hfil
 $\m@th\scriptstyle{##}$\hfil\crcr}
\def\endsmallmatrix{\crcr\egroup\egroup\,}
\newcount\c@MaxMatrixCols
\makeatother

\def\Z{\mathbb{Z}}
\def\Q{\mathbb{Q}}
\def\R{\mathbb{R}}
\def\C{\mathbb{C}}

\def\bs{\boldsymbol}

\title{Shintani correspondence for Maass forms of level $N$ 
and prehomogeneous zeta functions}
\author{Kazunari Sugiyama\!\!
\footnote{Department of Mathematics, Chiba Institute of Technology, 2-1-1 
Shibazono, Narashino, Chiba, 275-0023, Japan. 
E-mail:\texttt{skazu@sky.it-chiba.ac.jp}}
}
\date{\today}

\begin{document}

\maketitle

\begin{abstract}
A Shintani-Katok-Sarnak
type correspondence for Maass cusp forms of level $N$
is shown to be derived from analytic properties of
prehomogeneous zeta functions whose coefficients involve
periods of Maass  forms. 
\end{abstract}

\section*{Introduction}

In~\cite{Shimura73},  Shimura constructed a lifting from holomorphic cusp forms of
half-integral weight to cusp forms of integral weight. 
Shimura's original proof depends on the Rankin-Selberg method and
Weil's converse theorem~\cite{Weil}. 
In \cite{ShintaniLift}, Shintani constructed a lifting  from holomorphic cusp forms of 
integral weight to cusp forms of half-integral weight by using theta functions.
In the case of non-holomorphic modular forms, a prototype of
the lifting had already appeared in 
the work of Maa\ss~\cite{Maass}. 
Katok and Sarnak~\cite{KS} developed the method
of \cite{Maass} to prove 
 the Shintani correspondence for 
Maass cups forms of weight $0$ for $SL_2(\Z)$.
The Katok-Sarnak formula reveals a relation between 
the periods of Maass forms of weight $0$ and the 
Fourier coefficients of the corresponding form of weight $\frac12$, 
and now plays an important role in number theory.  
The Katok-Sarnak formula has been extended in many directions;  
we refer to Baruch-Mao~\cite{BaruchMao}, 
Biro~\cite{Biro}, 
Duke-Imamog\u{l}u-To\'{t}h~\cite{DIT}, 
Imamog\u{l}u-L\"{a}geler-To\'{t}h~\cite{ILT}. 
On the other hand, F.\ Sato~\cite{periodInd}
constructed a theory of prehomogeneous zeta functions
whose coefficients involve periods of automorphic forms.
In this note, we show that a Shintani-Katok-Sarnak
type correspondence is derived from analytic properties of a certain zeta 
functions 
investigated in~\cite{periodInd}. The proof relies on a Weil type converse theorem 
for Maass forms~\cite{MSSU}.

This is an announcement whose details will appear elsewhere. 

\section{Statement of the result}

The group $G= SL_{2}(\R)$ 
acts on the Poincar\'e upper half plane
$\mathcal{H}=\{z=x+iy\in \mathbb{C}\, |\, y>0\}$ via the
linear fractional transformation. 
Let $N$ be a positive integer and take a congruence subgroup 
$\Gamma_0(N)$ of level $N$ defined by
\[
\Gamma_0(N)=\{
\gamma\in SL_2(\Z) \, |\, \gamma_{21}\equiv 0 \pmod{N}\}, 
\]
where
$\gamma_{21}$ the $(2, 1)$-entry of $\gamma$.
Let $\chi$ be a Dirichlet character of mod $N$ satisfying 
$\chi(-1)=1$.  
We use the same symbol $\chi$ to denote the induced character of
$\Gamma_0(N)$ defined by $\chi(\gamma)= \chi(\gamma_{22})$ for 
$\gamma=(\gamma_{ij})\in \Gamma_0(N)$. 
A $C^{\infty}$-function $\Phi: \mathcal{H}\to \C$ is called  
a {\it Maass cusp form} of weight 0 for $\Gamma_0(N)$ with character $\chi$
if 
\begin{enumerate}
\item $\Delta_0 \Phi= \lambda(1-\lambda)\Phi$ for a $\lambda\in \C$, where
\[
 \Delta_0 = -y^2 \left(\frac{\partial^2}{\partial x^2}+
 \frac{\partial^2}{\partial y^2}\right)
\]
is the hyperbolic Laplacian on $\mathcal{H}$,  
\item $\Phi(\gamma z)=\chi(\gamma)\Phi(z)$ for $\gamma\in\Gamma_0(N)$, 
and
\item $F$ has exponential decay at all cusps of $\Gamma_0(N)$. 
\end{enumerate}
Let $\mathfrak{S}_0(N, \lambda, \chi)$ be the space of all such functions. 
For $g\in G$, we put $\phi(g)=\Phi(g^{-1}\cdot \sqrt{-1})$.
Let 
\[
 g= 
 \begin{pmatrix}
 \cos \theta & \sin \theta \\
 -\sin \theta & \cos \theta 
 \end{pmatrix}
 \begin{pmatrix}
 y^{-1/2} & 0 \\
 0 & y^{1/2}
 \end{pmatrix}
 \begin{pmatrix}
 1 & -x \\
 0 & 1
 \end{pmatrix}
\]
be the Iwasawa decomposition of $g\in G$. Then we have
$\phi(g)=\Phi(x+y\sqrt{-1})$. 
Let $V=\Sym_2(\R)$ be the space of real symmetric matrices of degree $2$.
Then $\widetilde{G}=\R^{\times}\times G$ acts on $V$ by 
$v \mapsto t \cdot g v \, {}^{t}\! g$ for $v\in V$ and $(t, g)\in \widetilde{G}$. 
Let $V_{+}=\{v\in V\, |\, \det v>0\}$ and 
$V_{-}=\{v\in V\, |\, \det v<0\}$. We have $V_{+}= \widetilde{G}\cdot I_2$
and $V_{-}=\widetilde{G}\cdot J_2$, where
\[
 I_2 = \begin{pmatrix}
 1 & 0 \\
 0 & 1
 \end{pmatrix}, \qquad 
 J_2 = \begin{pmatrix}
 0 & 1 \\
 1 & 0 
 \end{pmatrix}.
\]
Further, we put $H_{+}=SO(I_2)$ and $H_{-}=SO(J_2)$ so that
\begin{align*}
 H_{+} &=SO(2)=
 \left\{
 \left.
  k_{\theta}=
 \begin{pmatrix}
  \cos \theta & -\sin \theta \\
  \sin \theta & \cos \theta
  \end{pmatrix} \; \right|\; 
  0\leq \theta < 2\pi
 \right\}, \\
 H_{-} &=SO(1, 1) = 
 \left\{
  \left.
  a_{y}=
  \begin{pmatrix}
 y^{-1/2} & 0 \\
 0 & y^{1/2}
 \end{pmatrix}\; \right|\; 
 y\in \R^{\times}
 \right\}.
\end{align*}
We normalize the Haar measures $d\mu_{\pm}$ on $H_{\pm}$ by
\[
 d\mu_{+}(k_{\theta}) = \dfrac{d\theta}{2}, \qquad 
 d\mu_{-}( a_{y})= \dfrac{dy}{4|y|}.
\]
Let  $V_{+}^{p}$ (resp.\ $V_{+}^{n}$) be the set of positive (resp.\ negative)
definite symmetric matrices in $V_{+}$. For $v\in V$ with $\det v\neq 0$, 
we take $t_v >0$ and $g_v\in G$ such that  
\[
 v = \begin{cases}
 t_{v}(g_{v}I_2 \, {}^{t}\! g_{v}) & \text{if} \quad v\in V_{+}^{p}, \\
  -t_{v}(g_{v}I_2 \, {}^{t}\! g_{v}) & \text{if}\quad v\in V_{+}^{n}, \\
 t_{v}(g_{v}J_2 \, {}^{t}\! g_{v}) & \text{if}\quad v\in V_{-}. \\
 \end{cases}
\]
For $v\in V_{\Q}\cap V_{\pm}$, we define the {\it period}
$\mathcal{M}\phi(v)$ of $\phi$ by
\[
 \mathcal{M}\phi(v) = \int_{H_{\pm}/ g_{v}^{-1} \Gamma_{0, v} g_v}
 \phi(h g_{v}^{-1}) d\mu_{\pm}(h),
\]
where $\Gamma_{0, v}=\{\gamma \in \Gamma_0(N)\, |\, 
\gamma v \, {}^{t} \gamma = v\}$.
Then $\mathcal{M}\phi(v)$ is absolutely convergent and does not depend 
on the choice of $g_v$. 
By \cite[Lemma~6.3]{periodInd}, for  $v\in V_{+}$, we have
\[
 \mathcal{M}\phi(v) = \frac{\pi}{\varepsilon(v)} \cdot \Phi(z_{v}),
\]
where $\varepsilon(v)=\sharp (\Gamma_{0, v})$ and 
$z_{v}= g_{v}\cdot \sqrt{-1}$.
Note that $z_v$
coincides with  the so-called Heegner point associated with $v$. 
If $v\in V_{-}$, then $\{ g_v a_y\cdot \sqrt{-1}\, |\, y>0\}$ is 
the Heegner cycle associated with $v$, and thus 
$\mathcal{M}\phi(v)$ coincides (up to constant)
 with a certain cycle integral of $\Phi$.
Following the formulation of Shintani, we take a lattice 
$\mathcal{L}_{N}$ defined by
\[
 \mathcal{L}_{N}=
 \left\{ v=
 \begin{pmatrix}
 v_1 & N v_2 \\
 N v_2 & N v_3 
 \end{pmatrix} \; \bigg|\; 
 v_1, v_2, v_3\in \Z
 \right\}.
\]
(see \cite[p.\ 109]{ShintaniLift}). Further, for $v\in \mathcal{L}_{N}$, we put 
\[
 d_{N}(v):= N(v_2)^2 - v_1 v_{3} \quad (=-\frac{1}{N}\det v).
\]
Let $V_{\Z}$ be the set of half-integral  symmetric matrices of 
degree $2$: 
\[
 V_{\Z} = 
 \left\{ w^*=
 \begin{pmatrix}
 w_1^* & \frac{w_2^*}{2} \\
 \frac{w_{2}^*}{2} & w_3^* 
 \end{pmatrix} \; \bigg|\; 
 w_1^*, w_2^*, w_3^*\in \Z
 \right\},
\]
and for $w^*\in V_{\Z}$, we put 
\[
 \disc(w^*):= (w_2^*)^2 - 4 w_1^* w_{3}^* \quad (=-4 \det w^*).
\]
We take an automorphic factor $J(\gamma, z)$ of weight $\frac{1}{2}$ 
defined by
\[
 J(\gamma, z)=\frac{\theta(\gamma z)}{\theta(z)}, \quad
 \text{with} \quad  
\theta(z)=\sum_{n=-\infty}^{\infty} \exp(2\pi i n^2 z).
\]
Let
\[
  \Delta_{\frac12}
  =-y^2\left(\frac{\partial^2}{\partial x^2}+ \frac{\partial^2}{\partial y^2}\right)
  +\frac{iy}{2}\left(\frac{\partial}{\partial x}+ i\frac{\partial}{\partial y}\right)
\]
be the hyperbolic Laplacian of weight $\frac{1}{2}$ on $\mathcal{H}$, 
and $\psi$ a Dirichlet character of $\bmod 4N$.
A $C^{\infty}$-function $F: \mathcal{H}\to \C$ is called a 
Maass cusp form of weight $\frac{1}{2}$ for $\Gamma_0(4N)$ 
with character $\psi$ if
\begin{enumerate}
\item $\Delta_{\frac{1}{2}} F= \mu(1-\mu)F$ for a $\mu\in \C$, 
\item $F(\gamma z)=\psi(\gamma)J(\gamma, z)F(z)$ for $\gamma\in
\Gamma_0(4N)$, and
\item $F(z)$ has exponential decay at all cusps of $\Gamma_0(4N)$.
\end{enumerate}
We denote by $\mathfrak{S}_{\frac{1}{2}}(4N, \mu, \psi)$ the space of 
all such functions. Any $F\in  \mathfrak{S}_{\frac{1}{2}}(4N, \mu, \psi)$
has a Fourier expansion of the form 
\begin{equation}
 \label{form:FourierExp-wt1/2}
 F(z)= \sum
 \begin{Sb}
 n=-\infty \\
 n\neq 0
 \end{Sb}^{\infty} c(n) \cdot W_{1, \mu}(n, y) 
 \boldsymbol{e}[nx],
\end{equation}
where $\boldsymbol{e}[x]= \exp(2\pi \sqrt{-1} x)$ and
for $\ell\in \Z$, 
\begin{equation}
 \label{form:Whittaker}
W_{\ell, \mu}(n, y)=
y^{-\frac{\ell}{4}} W_{\frac{\sgn(n)\ell}{4}, \, \mu-\frac{1}{2}}
 \left(4\pi|n|y\right).
\end{equation}
Here 
$W_{\kappa, \nu}(z)$ denotes the Whittaker function. 
For a Dirichlet character $\chi$ of $\bmod{N}$, let 
\[
 \tau_{\chi}(n)= \sum
 \begin{Sb}
 m\bmod{N} \\
 (m, N)=1
 \end{Sb} \chi(m) \boldsymbol{e}\left[\frac{mn}{N}\right]
\]
be the Gauss sum associated with $\chi$.
Now we state our main theorem. 
\begin{thm}
 \label{thm:MainThm}
Let $\lambda\neq \frac{1}{2}$ and 
assume that 
$\Phi(z)\in \mathfrak{S}_{0}(N, \lambda, \chi^2)$. We put
\[
\mu=\frac{2\lambda+1}{4}, \qquad
\chi_N(r) = \chi(r) \left(\dfrac{N}{r}\right).
\]
Then there exists an $F(z)\in \mathfrak{S}_{\frac12}\left(4N, \mu, \chi_{N}\right)$
such that the Fourier coefficients $c(n)$ in \eqref{form:FourierExp-wt1/2}
are given by
\begin{align*}
 c(n) &=
 2 \pi^{-\frac12} \cdot n^{-\frac34} 
\sum
\begin{Sb}
 v\in \Gamma_0(N) \backslash \mathcal{L}_{N} \\
 d_{N}(v)=n
\end{Sb}
{\chi(v_1)}
\mathcal{M}\phi(v), \\ 
c(-n) &=
n^{-\frac34}
\sum
\begin{Sb}
 v\in \Gamma_0(N) \backslash \mathcal{L}_{N}\\
d_{N} (v)=-n
\end{Sb}
\frac{{\chi(v_1)}\Phi(z_v)}{\epsilon(v)},
 \end{align*}
 for $n=1, 2, 3,\dots$.
 Furthermore, if we put
\begin{align*}
&c^*(n) =
2^{\lambda}\pi^{-\frac{1}{2}}\cdot 
n^{-\frac34}
\sum
\begin{Sb}
 w^*\in \Gamma_0(N) \backslash V_{\Z}  \\
 \disc w^*=n
\end{Sb}
{\tau_{{\chi}}(w_3^*)\mathcal{M}\phi(w^*)}, \\
&c^{*}(-n) = 2^{\lambda-1}\cdot 
n^{-\frac34}
\sum
\begin{Sb}
 w^*\in \Gamma_0(N) \backslash V_{\Z}  \\
 \disc w^*=-n
\end{Sb}
\frac{\tau_{{\chi}}(w_3^*)\Phi(z_{w^*})}
{\epsilon(w^*)},
\end{align*}
for $n=1,2,3,\dots$, 
and define a function $G(z)$ on $\mathcal{H}$ by
\[
G(z)= N^{-\frac{3}{4}}\sum
 \begin{Sb}
 n=-\infty \\
 n\neq 0
 \end{Sb}^{\infty} c^*(n) \cdot 
 W_{1, \mu}(n, y) \boldsymbol{e}[nx], 
\]
then we have 
$G(z)\in \mathfrak{S}_{\frac12}\left(4N, \mu, \overline{\chi}\right)$
and 
\[
F\left(-\frac{1}{4Nz}\right)(\sqrt{N}z)^{-\frac12} =
\boldsymbol{e}\left[-\frac18\right]
\cdot G(z).
\]
\end{thm}

\section{A Weil type converse theorem for Maass cusp forms}

Our proof for Theorem~\ref{thm:MainThm} relies on a converse 
theorem given in~\cite{MSSU}. 
Here let us recall briefly the result, with some modifications. 
For the convenience of readers, we give the statement for 
general weights.
Fix an integer $\ell$ and a positive integer  $N$. 
We assume that $N$ is a multiple of $4$ when $\ell$ is odd.
Let 
 $\alpha=\{\alpha(n)\}_{n\in \Z\setminus\{0\}}$ and $\beta=\{\beta(n)\}_{n\in \Z\setminus\{0\}}$ be  complex sequences of polynomial growth. 
For $\alpha,\beta$, we can define the  $L$-functions 
 $\xi_{\pm}(\alpha;s), \xi_{\pm}(\beta;s)$ by
\[
 \xi_{\pm}(\alpha;s) = 
 \sum_{n=1}^{\infty}\frac{\alpha(\pm n)}{n^s}, \quad
  \xi_{\pm}(\beta;s) = 
 \sum_{n=1}^{\infty}\frac{\beta(\pm n)}{n^s},
\]
and the completed $L$-functions 
$\Xi_{\pm}(\alpha;s)$ and $\Xi_{\pm}(\beta;s)$
by
$\Xi_{\pm}(\alpha;s)=
 (2\pi)^{-s} \Gamma(s) \xi_{\pm}(\alpha;s)$ and
$\Xi_{\pm}(\beta;s)=
 (2\pi)^{-s} \Gamma(s) \xi_{\pm}(\beta;s)$.

Now we assume the following two conditions: 

\begin{description}
 \item[{\bf [C1]}] The $L$-functions $\xi_{\pm}(\alpha;s), \xi_{\pm}(\beta;s)$  have analytic continuations to {\it entire} functions of $s$, and 
 are of finite order in any vertical strip. 
\end{description}

\begin{description}
 \item[{\bf [C2]}] The following functional equation holds:
 \begin{multline}
  \label{form:AssumpFE}
  \gamma(s)\left(
 \begin{array}{c}
 \Xi_{+}(\alpha;s) \\[2pt] \Xi_{-}(\alpha;s)
 \end{array}
 \right) 
 = N^{2-2\mu-s}\cdot 
 \Sigma(\ell) 
 \cdot
 \gamma(2-2\mu-s)\left(
 \begin{array}{c}
 \Xi_{+}(\beta;2-2\mu-s) \\[2pt] \Xi_{-}(\beta;2-2\mu-s)
 \end{array}
 \right),
 \end{multline}
 where $\gamma(s)$ and $\Sigma(\ell)$ are given by
 \[
  \gamma(s) =
  \begin{pmatrix}
  e^{\pi si/2} & e^{-\pi s i/2} \\
  e^{-\pi si/2} & e^{\pi si/2}
  \end{pmatrix}, \quad
  \Sigma(\ell)=
  \begin{pmatrix}
  0 & i^{\ell} \\
  1 & 0
  \end{pmatrix}.
 \]
\end{description}


For an odd prime number $r$ with $(N,r)=1$ and a Dirichlet character $\psi$ mod $r$, the twisted $L$-functions 
  $\xi_{\pm}(\alpha, \tau_{\psi};s), 
    \xi_{\pm}(\beta, \tau_{\psi};s)$ 
are defined by
 \begin{align*}
  \xi_{\pm}(\alpha, \tau_{\psi};s) &= 
 \sum_{n=1}^{\infty}\frac{\alpha(\pm n)\tau_{\psi}(\pm n)}{n^s}, \\
  \xi_{\pm}(\beta,\tau_{\psi};s) &= 
 \sum_{n=1}^{\infty}\frac{\beta(\pm n)\tau_{\psi}(\pm n)}{n^s}, 
 \end{align*}
where $\tau_{\psi}(n)$ is the Gauss sum associated with $\psi$. 
The complete $L$-functions 
$\Xi_{\pm}(\alpha,\tau_{\psi}; s)$ and
$\Xi_{\pm}(\beta, \tau_{\psi};s)$ are defined by
$\Xi_{\pm}(\alpha,\tau_{\psi}; s)=
 (2\pi)^{-s} \Gamma(s) \xi_{\pm}(\alpha,\tau_{\psi};s)$ and
$\Xi_{\pm}(\beta, \tau_{\psi};s)=
 (2\pi)^{-s} \Gamma(s) \xi_{\pm}(\beta,\tau_{\psi};s)$, respectively.

Let $\mathbb{P}_N$ be a set of odd prime numbers not dividing $N$ such that, for any positive integers $a, b$ coprime to each other, $\mathbb{P}_N$ contains a prime number $r$ of the form $r=am+b$ for some $m \in \Z_{>0}$.   
For an $r \in \mathbb{P}_N$, denote by $X_r$ the set of all Dirichlet characters mod $r$ (including the principal character). 
For $\psi \in X_r$, we define the Dirichlet character  $\psi^*$  by
\begin{equation}
\psi^*(k)=\overline{\psi(k)}\left(\frac{k}{r}\right)^{\ell}. 
\label{eqn:psi-star}
\end{equation}
For an odd integer $d$, we put $\varepsilon_d =1$ or $\sqrt{-1}$ according as
$d\equiv 1$ or $3 \pmod{4}$. 
Let
\[
 C_{\ell,r} =
  \begin{cases}
  1 & (\text{$\ell$ is even}), \\
  \varepsilon_{r}^{\ell}  &   (\text{$\ell$ is odd}).
 \end{cases}
\]

In the following, we fix a Dirichlet character $\chi$  mod ${N}$
 that satisfies $\chi(-1)=(\sqrt{-1})^{\ell}$ (resp.\ $\chi(-1)=1$)
 when $\ell$ is even (resp.\ odd).

For an $r \in \mathbb{P}_N$ and a $\psi \in X_r$, we consider the following conditions $\mathrm{[C1]}_{r,\psi}$ -- $\mathrm{[C2]}_{r,\psi}$ on  
$\xi_{\pm}(\alpha,\tau_{\psi};s)$ and $\xi_{\pm}(\beta,\tau_{\psi^*};s)$.  
 
\begin{description}
\item[$\text{\bf [C1]}_{r,\psi}$] 
 $\xi_{\pm}(\alpha,\tau_{\psi};s), \xi_{\pm}(\beta,\tau_{\psi^*};s)$  have analytic 
 continuations to {entire} functions of $s$, and are of finite order in any vertical strip. 
\item[$\text{\bf [C2]}_{r,\psi}$] 
$\Xi_{\pm}(\alpha,\tau_{\psi};s)$ and $\Xi_{\pm}(\beta,\tau_{\psi^*};s)$
satisfy the following functional equation:
\begin{multline}
\label{form:AssumpTwistedFE}
 \gamma(s)\left(
 \begin{array}{c}
 \Xi_{+}(\alpha,\tau_{\psi};s) \\ \Xi_{-}(\alpha,\tau_{\psi};s)
 \end{array}
 \right) 
  = {\chi(r)}\cdot C_{\ell,r}
   \cdot  \psi^*(-N) 
  \cdot 
 r^{2\mu-2}\cdot (Nr^2)^{2-2\mu-s} 
 \cdot \Sigma(\ell) \\
 \cdot  \gamma(2-2\mu-s)
 \begin{pmatrix}
 \Xi_{+}\left(\beta, \tau_{{\psi}^*};2-2\mu-s\right) \\[4pt]
 \Xi_{-}\left(\beta, \tau_{{\psi}^*};2-2\mu-s\right)
 \end{pmatrix}. 
 \end{multline}
\end{description}
 
\begin{lem}
 \label{corollary:Maassforms}
Let  $\mu\not\in \frac12\Z$. 
We assume that $\xi_\pm(\alpha;s)$ and $\xi_\pm(\beta;s)$ satisfy the  conditions
 {\rm [C1]} and {\rm  [C2]}.
We assume furthermore that, for any $r \in \mathbb{P}_N$ and $\psi \in X_r$, 
 $\xi_{\pm}(\alpha,\psi;s)$ and $\xi_{\pm}(\beta,\psi^*;s)$  satisfy the conditions $\mathrm{[C1]}_{r,\psi}$ and $\mathrm{[C2]}_{r,\psi}$.
 We define the function $\widetilde{W}_{\ell, \mu}(n, y)$ by 
 \[
 \widetilde{W}_{\ell, \mu}(n, y)=
 \frac{|n|^{\mu-1}}
{\Gamma\left(\mu+\frac{\sgn(n)\ell}{4}\right)}  
\cdot W_{\ell, \mu}(n, y),
\]
where $W_{\ell, \mu}(n, y)$ is given as \eqref{form:Whittaker}, 
and the functions $F_\alpha(z)$ and  $G_\beta(z)$ on $\mathcal{H}$ by 
\begin{align*}
 F_\alpha(z) &
  =\sum
 \begin{Sb}
 n=-\infty \\
 n\neq 0
 \end{Sb}^{\infty} \alpha(n)\cdot  
 \widetilde{W}_{\ell, \mu}(n, y)
 \mathbf{e}[nx],  \\
G_\beta(z) &
  = N^{1-\mu}  \sum
 \begin{Sb}
 n=-\infty \\
 n\neq 0
 \end{Sb}^{\infty} \beta(n) \cdot 
  \widetilde{W}_{\ell, \mu}(n, y)
 \mathbf{e}[nx]. 
\end{align*}
%
%
Then $F_\alpha(z)$ (resp.\ $G_\beta(z)$) gives a Maass cusp
form 
 for ${\Gamma}_0(N)$ of weight $\frac{\ell}{2}$ with
character $\chi$ (resp.\ $\chi_{N,\ell}$), 
and eigenvalue $(\mu-\ell/4)(1-\mu-\ell/4)$, where
\[
\chi_{N,\ell}(d)=\overline{\chi(d)}\left(\frac{N}{d}\right)^\ell.
\]
Moreover, 
we have 
\[
F_{\alpha}\left(-\frac{1}{Nz}\right)(\sqrt{N}z)^{-\ell/2} =G_{\beta}(z).
\]
\end{lem}

\begin{remark}
Here we have assumed a stronger condition 
$\mu\not\in \frac12 \Z$ than that given in 
the previous paper~\cite{MSSU}. 
This enables us to remove conditions on
zeros of $L$-functions (cf.\  \cite[p.\ 33]{MSSU}).  
We also note that the cuspidality of
$F_{\alpha}(z), G_{\beta}(z)$ follows from 
the entireness of $L$-functions.
\end{remark}

\section{Prehomogeneous zeta functions}

As an example of the theory of~\cite{periodInd},  Sato 
investigated the zeta functions associated to the vector space
of symmetric matrices of degree 2 whose coefficients involve
the periods $\mathcal{M}\phi(v)$ of Maass cusp forms $\Phi$. 
In this section, we introduce twisted versions of these zeta functions
and give their analytic properties such as analytic continuations and 
functional equations.

Keep the notation as in the previous sections. We define  zeta functions
$\zeta_{\pm}(\phi, \chi;s)$ and
 $\zeta_{\pm}^{*}(\phi, \tau_{\chi};s)$ by
\begin{align*}
 \zeta_{\pm}(\phi, \chi;s) &=
 \sum
  \begin{Sb}
  v\in \Gamma_0(N)\backslash \mathcal{L}_{N} \\
  \sgn d_{N}(v)= \pm
  \end{Sb}
   \frac{\chi(v_1) \mathcal{M}\phi(v)}{|d_{N}(v)|^s}, \\
  \zeta_{\pm}^{*}(\phi, \tau_{\chi};s) &=
 \sum
  \begin{Sb}
  w^*\in \Gamma_0(N)\backslash V_{\Z} \\
  \sgn \disc (w^*)= \pm
  \end{Sb}
   \frac{\tau_{\chi}(w_3^*) \mathcal{M}\phi(w^*)}{|\disc w^*|^s}.
\end{align*}
Then we have the following lemma, whose proof is similar to that of
\cite[Theorem~6.7]{periodInd}.
\begin{lem}
The zeta functions $\zeta_{\pm}(\phi, \chi;s)$ and
$\zeta_{\pm}^{*}(\phi, \tau_{\chi};s)$ have analytic continuations to 
entire functions of $s$ and satisfy the following functional
equation:
\begin{multline}
\label{form:FEofPVLfunction}
\begin{pmatrix}
\zeta_{+}(\phi, \chi; \frac{3}{2}-s) \\[3pt]
\zeta_{-}(\phi, \chi; \frac{3}{2}-s)
\end{pmatrix} 
=\pi^{\frac12-2s}  N^{s-\frac32}\, \Gamma
\left(s+\frac{\lambda-1}{2}\right) \Gamma
\left(s-\frac{\lambda}{2}\right) 
\cdot \Psi_{\lambda}(s)
\begin{pmatrix}
\zeta_{+}^*(\phi, \tau_{\chi}; s) \\[3pt]
\zeta_{-}^*(\phi, \tau_{\chi}; s)
\end{pmatrix},
\end{multline}
where $\Psi_{\lambda}(s)$ is a $2\times 2$ matrix given by
\[
\Psi_{\lambda}(s)=
\begin{pmatrix}
\sin \pi s &
\displaystyle
\frac{2^{\lambda-1} \cdot \pi\Gamma(1-\lambda)}
{ \Gamma\left(1-\frac{\lambda}{2}\right)^2}
\cos\frac{\pi \lambda}{2} \\
\displaystyle
\frac{\Gamma\left(1-\frac{\lambda}{2}\right)^2}{2^{\lambda-1} \cdot 
\pi\Gamma(1-\lambda)}
\sin\frac{\pi \lambda}{2} &
\cos\pi s
\end{pmatrix}.
\]

\end{lem}

Let $r$ be an odd prime number $r$ with $(N,r)=1$, and 
$\psi$ a Dirichlet character of mod $r$. We denote by $\psi^*$ 
the Dirichlet character defined  as
\eqref{eqn:psi-star} with $\ell =1$.
We define 
$\zeta_{\pm}(\phi, \chi, \tau_{\psi};s)$ and
$\zeta_{\pm}^{*}(\phi, \tau_{\chi}, \tau_{\psi^*};s)$ by
\begin{align*}
 \zeta_{\pm}(\phi, \chi, \tau_{\psi};s) 
 &=
 \sum
  \begin{Sb}
  v\in \Gamma_0(N)\backslash \mathcal{L}_{N} \\
  \sgn d_{N}(v)= \pm
  \end{Sb}
   \frac{\chi(v_1) 
   \mathcal{M}\phi(v) \tau_{\psi}(d_N(v))}{|d_{N}(v)|^s}, \\
  \zeta_{\pm}^{*}(\phi, \tau_{\chi}, \tau_{\psi^*};s) 
  &=
 \sum
  \begin{Sb}
  w^*\in \Gamma_0(N)\backslash V_{\Z} \\
  \sgn \disc (w^*)= \pm
  \end{Sb}
   \frac{\tau_{\chi}(w_3^*)
    \mathcal{M}\phi(w^*)
    \tau_{\psi^*}(\disc w^*)}{|\disc w^*|^s}.
\end{align*}
Then we have the following lemma.
\begin{lem}
\label{lem:ZetaTwistedbyGS}
The zeta functions $\zeta_{\pm}(\phi, \chi, \tau_{\psi};s)$ and
$\zeta_{\pm}^{*}(\phi, \tau_{\chi}, \tau_{\psi^*};s)$ have analytic continuations to 
entire functions of $s$ and satisfy the following functional
equation:
\begin{multline}
\label{form:FEoftwistedPV}
\begin{pmatrix}
\zeta_{+}(\phi, \chi, \tau_{\psi}; \frac{3}{2}-s) \\[3pt]
\zeta_{-}(\phi, \chi, \tau_{\psi}; \frac{3}{2}-s)
\end{pmatrix} 
= \varepsilon_{r} \chi_{N}(r) \psi^*(-4N) 
\cdot r^{2s-\frac32}
\pi^{\frac12-2s}  N^{s-\frac32}
\, \Gamma
\left(s+\frac{\lambda-1}{2}\right) \Gamma
\left(s-\frac{\lambda}{2}\right) \\
\cdot 
\Psi_{\lambda}(s)
\begin{pmatrix}
\zeta_{+}^*(\phi, \tau_{\chi},\tau_{\psi^*}; s) \\[3pt]
\zeta_{-}^*(\phi, \tau_{\chi}, \tau_{\psi^*}; s)
\end{pmatrix}.
\end{multline}
\end{lem}

The proof of Lemma~\ref{lem:ZetaTwistedbyGS} goes along
the same line as Sato~\cite{ZetaDist}, Ueno~\cite{Ueno}.
In this case, however, it is necessary to calculate a kind of 
Gauss sums that have not appeared in the previous works. 
The author has learned such calculation from unpublished notes
of Sato.  
We quote his result, which is a key ingredient and of independent interest. 
Let $f_{\psi, \chi}(v)$ be a function on $V_{\Q}$ defined by
\[
 f_{\psi, \chi}(v) = 
 \begin{cases}
 \tau_{\psi}(d_N(v)) \cdot \chi(v_1) & (v\in \mathcal{L}_{N}) \\
 0 & (v\not\in \mathcal{L}_{N})
 \end{cases}.
\]
Let $\langle v, v^*\rangle$ be the inner product on $V$ defined
by $\langle v, v^*\rangle=\tr(v w v^* w^{-1})$ with
$w=
\left(\begin{smallmatrix}
0 & 1\\
-1 & 0
\end{smallmatrix}\right)$.
For $v^*\in V_{\Q}$, 
we define the Fourier transform $\widehat{f_{\psi, \chi}}(v^*)$ by
\begin{equation}
\label{form:DefOfFTofFpsichi}
\widehat{f_{\psi, \chi}}(v^*) =
\frac{1}{[V_{\Z}:L]} \sum_{v\in V_{\Q}/L}
f_{\psi, \chi}(v) \bs{e}[\langle v, v^*\rangle], 
\end{equation}
where $L$ is a sufficiently small lattice so that $L\subset V_{\Z}$ and
the value $f_{\psi, \chi}(v) \bs{e}[\langle v, v^*\rangle]$ depends only on
the residue class $v+L$. 
\begin{lem}[F.\ Sato]
If $v^*\not\in \frac{1}{Nr}V_{\Z}$, then we have 
$\widehat{f_{\psi, \chi}}(v^*) =0$. 
If $v^* =\frac{1}{Nr} w^*\in \frac{1}{Nr}V_{\Z}$, we have
\[
\widehat{f_{\psi, \chi}}(v^*) 
=\frac{\varepsilon_r}{2r^{\frac{3}{2}}N^3} 
 \chi_N(r) \cdot
 \psi^{*}(-4N) 
 \tau_{{\chi}}(w_3^*) \tau_{\psi^*}(\disc(w^*)).
\]
\end{lem}

\section{An outline of the proof of Theorem~\ref{thm:MainThm}}

We construct $L$-functions satisfying 
two conditions [C1] and [C2]. 
In the functional equation~\eqref{form:AssumpFE},
we let $\ell=1$ and $\mu= \frac{2\lambda+1}{4}$, and
replace $N$ by $4N$. 
Then it follows from an elementary calculation that 
\eqref{form:AssumpFE} is transformed as
\begin{multline}
\label{form:transformed desired FE of xi}
\begin{pmatrix}
\xi_{+}(\alpha; s) \\
\xi_{-}(\alpha; s) 
\end{pmatrix} =(4N)^{\frac{3}{2}-\lambda-s}\cdot 
2^{2s+\lambda-\frac{3}{2}}
\cdot \pi^{2s+\lambda-\frac{5}{2}}  
\cdot 
\boldsymbol{e}\left[\frac18\right]  
\Gamma(1-s)\Gamma\left(\frac{3}{2}-\lambda-s\right) \\ 
\cdot
\begin{pmatrix}
-\cos \pi (s+\frac{\lambda}{2}) & 
\sin \frac{\pi\lambda}{2}\\
 \cos\frac{\pi \lambda}{2} &
-
\sin\pi (s+\frac{\lambda}{2})
\end{pmatrix} 
\cdot 
\begin{pmatrix}
\xi_{+}(\beta; \frac{3}{2}-\lambda-s) \\[4pt]
\xi_{-}(\beta; \frac{3}{2}-\lambda-s) 
\end{pmatrix}.
\end{multline} 
We put
\begin{align*}
\widetilde{\zeta}_{+}(\phi, \chi; s) &:= 2^{2-\lambda}\cdot 
\frac{\Gamma(\lambda)}{\Gamma(\frac{\lambda}{2})^2}\cdot 
\zeta_{+}\left(\phi, \chi; s+\frac{\lambda}{2}\right), \\
\widetilde{\zeta}_{-}(\phi, \chi; s) &:=
\zeta_{-}\left(\phi, \chi; s+\frac{\lambda}{2}\right), \\
\widetilde{\zeta}_{+}^*(\phi, \tau_{\chi}; s) &:= 2^{\frac12}\cdot
N^{-\frac{3}{2}+\frac{\lambda}{2}}\cdot  
\frac{\Gamma(\lambda)}{\Gamma(\frac{\lambda}{2})^2} \cdot 
\zeta_{+}^*\left(\phi, \tau_{\chi}; s+\frac{\lambda}{2}\right), \\
\widetilde{\zeta}_{-}^*(\phi, \tau_{\chi}; s) &:= 
2^{\lambda-\frac{3}{2}}\cdot 
N^{-\frac{3}{2}+\frac{\lambda}{2}}\cdot  
\zeta_{-}^*\left(\phi, \tau_{\chi}; s+\frac{\lambda}{2}\right).
\end{align*}
Then \eqref{form:FEofPVLfunction}
can be rewritten as
\begin{multline}
\begin{pmatrix}
\widetilde{\zeta}_{+}(\phi, \chi; s) \\
\widetilde{\zeta}_{+}(\phi, \chi; s) 
\end{pmatrix} =(4N)^{\frac{3}{2}-\lambda-s}\cdot 
2^{2s+\lambda-\frac{3}{2}} 
\cdot \pi^{2s+\lambda-\frac{5}{2}}  
\cdot 
\boldsymbol{e}\left[\frac18\right]  
\Gamma(1-s)\Gamma\left(\frac{3}{2}-\lambda-s\right)  \\
\cdot
\begin{pmatrix}
-\cos \pi (s+\frac{\lambda}{2}) & 
\sin \frac{\pi\lambda}{2}\\
 \cos\frac{\pi \lambda}{2} &
-
\sin\pi (s+\frac{\lambda}{2})
\end{pmatrix} 
\cdot 
\begin{pmatrix}
\widetilde{\zeta}_{+}^*(\phi, \tau_{\chi}; \frac{3}{2}-\lambda-s) \\[4pt]
\widetilde{\zeta}_{-}^*(\phi, \tau_{\chi}; \frac{3}{2}-\lambda-s) 
\end{pmatrix},
\end{multline}
which agrees with~\eqref{form:transformed desired FE of xi}.
Similarly, the functional equation 
\eqref{form:FEoftwistedPV} of the twisted zeta functions 
can be compared to~\eqref{form:AssumpTwistedFE} in 
the condition [C2]${}_{r, \psi}$.
Now the converse theorem (Lemma~\ref{corollary:Maassforms}) applies, 
and Theorem~\ref{thm:MainThm} is obtained. 
Further details will be discussed elsewhere.

\begin{remark}
In a paper~\cite{ILT} that appeared very recently, the Katok-Sarnak formula 
is generalized for Maass forms of even weight and odd level 
with trivial characters.
It is an interesting problem to combine their technique,
such as use of differential operators, with our method.  
\end{remark}

\subsection*{Acknowledgement.} The author wishes to thank
Professor Fumihiro Sato for his kind guidance and 
for giving the author permission to use his unpublished results.

\end{document}